# THE SPIRAL VORTEX

# A story about tornados and bathtubs

Khristo N. Boyadzhiev
Department of Mathematics and Statistics
Ohio Northern University
Ada, Ohio 45810
k-boyadzhiev@onu.edu

## Introduction

When we drain the water from the bathtub we see a small whirlpool. Amazingly, some galaxies and hurricanes look like a whirlpool too (Figures 1 and 2).

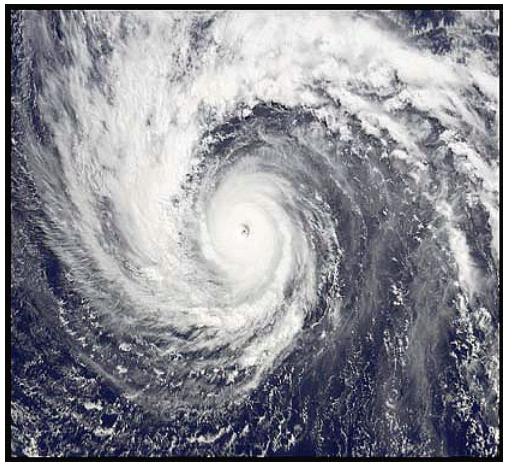

**Figure 1 Hurricane Bonnie, August 1998**

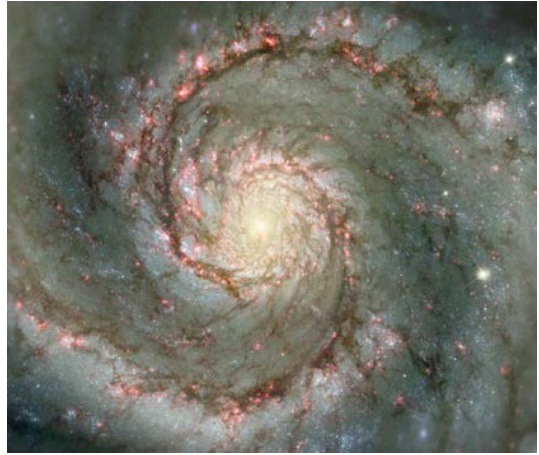

**Figure 2 Spiral Galaxy M51**

This likeness is not coincidental - the spiral lines observed in such cases closely match equiangular spirals (Figure 3).

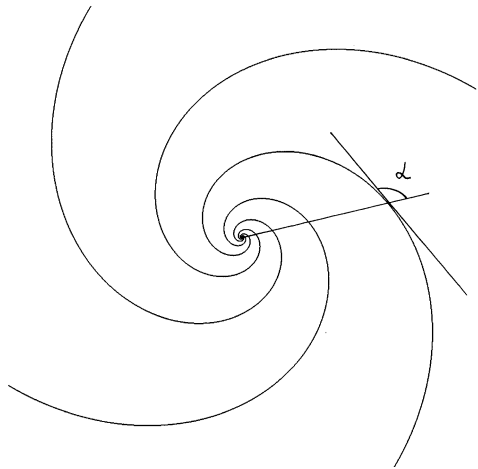

**Figure 3 Equiangular spiral**

The popular equiangular spiral is defined by the polar equation

$$r(\theta) = r_0 e^{k\theta} \ (-\infty < \theta < \infty)$$

It has the characteristic property that the angle $\alpha$ between the radius vector and the tangent vector is the same at any point on the curve. In the above equation, $k = \cot\alpha$. This spiral appears in many formations in nature [2, 5, 7, 8, 13]. On Figure 3 we see several spirals with the same angle $\alpha$ and different values of $r_0$.

Using complex variables, we shall give a simple explanation why vortex curves are equiangular spirals and not something else. First we consider two-dimensional velocity vector fields with one singularity at the origin $O$. Such fields can model a water surface, or a thin layer in the atmosphere. We also involve systems of differential equations which help to extend our simple model from two to three dimensions.

This article has an educational character and is designed for undergraduate student with background in complex variables and differential equation. The author believes it will be interesting also to all nature-curious mathematicians.

## Vector fields on the plane

For convenience, we identify the $xy$-plane with the complex $z$-plane by setting $z = x + iy = (x, y)$. Vectors in standard position (starting from the origin), and complex numbers are also identified, as in equation (1) below.

Suppose now we have a fluid on the complex plane, described by the velocity vector field

$$V(z) = V(x, y) = <u, v> = u + iv, \qquad (1)$$

defined everywhere except, possibly, at the origin. Every point (particle) $M(x, y)$ has velocity $V = <u(x, y), v(x, y)>$. Governed by this velocity vector field, the particles move on streamlines. A streamline (trajectory) is a curve, at each point of which, the velocity vector is tangent to that curve. If the smooth curve $L$ is defined by the parametric equations

$$x = x(t), \ y = y(t), \ t_1 \leq t \leq t_2, \qquad (2)$$

then $L$ is a streamline for the velocity vector field (1) if and only if the tangent vector $<x', y'>$

and the velocity vector $<u, v>$ are parallel at every point $(x, y)$ on $L$.

Some simple definitions from vector calculus are needed. Let $G$ be a closed, positively oriented (counterclockwise) smooth curve surrounding a convex domain $D$. We can think that $D$ is a just a disk. The integral

$$Cir(V:G) = \oint_G V dr = \oint_G u dx + v dy \qquad (3)$$

(where $dr = <dx, dy>$) defines the *circulation* of the vector field along the curve. Nonzero circulation indicates the presence of whirls inside $D$.

We also define the *flux* through the curve $G$

$$Flux(V:G) = \oint_G u dy - v dx. \qquad (4)$$

Nonzero flux indicates the appearance or disappearance of fluid inside $D$, i.e. the presence of sources or sinks. We now formally define a *vortex,* a *source* and a *sink*.

**Definition 1.** A point $M$ is called a *vortex* for the vector field $V$, if there is a neighborhood $U$ of $M$ such that the circulation $Cir(V:G)$ on any circle $G \subseteq U$ centered at $M$ is nonzero.

Shortly, at a vortex the fluid spins, not necessarily appearing or disappearing.

**Definition 2**. With $M$ and $U$ as above, the point $M$ is called a *source*, if $Flux(V:G) > 0$ and a *sink*, if $Flux(V:G) < 0$ for every circle $G \subseteq U$ centered at $M$.

In these definitions we do not require $M$ to be in the domain of the field. It may be a singular point.

Before proceeding further, recall that multiplying one complex number $z$ by the exponential $e^{i\alpha}$ produces a counterclockwise rotation about $O$ by the angle with radian measure $\alpha$.

When considering a flat fluid with possibly one source/sink, it is reasonable to assume that the fluid is described by a two-dimensional velocity vector field $V$ with possibly one singularity at the origin, which may be a sink, or a source, and/or a vortex (or neither), and there are no other sinks or sources except $O$. We also do not permit the presence of any whirls which are not centered at $O$. It is natural to assume that this field is rotationally invariant with respect to the origin $O$. Here are the exact conditions:

(A)   The vector field $V = u + iv$ is smooth, i.e. $u, v$ have continuous partial derivatives everywhere except, possibly, at the origin $O$.

(B)   The vector field has no sinks or sources anywhere except, possibly, at the origin $O$, and also no whirls. This means both integrals (3), (4) are zero for every closed curve whose interior is separated from $O$.

(C)   $V$ is rotationally invariant, that is, $e^{i\alpha}V(z) = V(e^{i\alpha}z)$ for all $z \neq 0$, $\alpha \in (0, \pi)$. When we rotate the plane about the origin at angle $\alpha$, the picture we see does not change.

Now we will give a simple characterization of such vector fields and will show that their streamlines are equiangular spirals.

**Spiral Vortex Theorem.** A vector field $V$ satisfies (A), (B) and (C), if and only if it has the form

$$V(z) = c / \overline{z} \tag{5}$$

where $c = a + ib$ is a complex constant. When $a \neq 0$ and $b \neq 0$, the streamlines of this field are equiangular spirals with polar equation

$$r = r_0 e^{kt}, -\infty < t < \infty \tag{6}$$

where $r_0 = r(0)$ is an arbitrary real constant (initial condition) and $k = a/b$. When $b = 0$ the streamlines are rays going from $O$ to infinity (Figure 5 below). When $a = 0$ the stream lines are concentric circles centered at the origin (Figure 4).

*Proof of the theorem.* Suppose the vector field $V(z) = u + iv$ (where $u = u(x, y), v = v(x, y)$), satisfies the above conditions. Let $D$ be an arbitrary disk separated from the origin with boundary $G$. Then condition (B) implies in view of Green's theorem.

$$Cir(V : G) = \oint_G u\,dx + v\,dy = \iint_D (v_x - u_y)\,dx\,dy = 0$$

and also

$$\text{Flux}(V:G) = \oint_G u\,dy - v\,dx = \oiint_D (u_x + v_y)\,dx\,dy = 0.$$

Since the disk $D$ is arbitrary, we conclude that

$$u_x - v_y = 0, \quad u_y + v_x = 0$$

everywhere except possibly at the origin. These two equations are the Cauchy-Riemann equations for the function $\overline{V}(z) = u - iv$. Therefore, this function is holomorphic on the entire complex plane indented at the origin. It has a Laurent series convergent for every $z \neq 0$,

$$\overline{V}(z) = \sum_{n=-\infty}^{\infty} c_n z^n.$$

According to property (C), we have $e^{-i\alpha}\overline{V}(z) = \overline{V(e^{i\alpha}z)}$ for every $\alpha \in (0, \pi)$. This gives $e^{-i\alpha}c_n = e^{i\alpha n}c_n$, or $c_n = e^{i\alpha(n+1)}c_n$ for every integer $n$. This implies $c_n = 0$ for every $n \neq -1$ and therefore, $\overline{V(z)} = c_{-1}/z$. We conclude that

$$V(z) = c/\overline{z}$$

for the complex constant $c = \overline{c}_{-1}$.

One surprising result from this theorem is that the vector field $V$ is zero at infinity!

It is clear from the above considerations that any vector field of the form (5) satisfies (A), (B), and (C).

We continue now with the proof to show that the streamlines of (5) are equiangular spirals. Let $c = a + ib$ where $a, b$ are real. With $z = x + iy$ equation (5) can be written in the form:

$$V(z) = V(x.y) = \left\langle \frac{ax - by}{|z|^2}, \frac{bx + ay}{|z|^2} \right\rangle.$$

Consider the system of differential equations:

$$\begin{aligned} x' &= ax - by \\ y' &= bx + ay \end{aligned} \tag{7}$$

where $x = x(t), y = y(t)$. The vector field $<x', y'>$ defined by this system is parallel to $V(x, y)$ at each point, therefore the streamlines of $V$ are determined by the solutions of (7). When $b \neq 0$ this system is equivalent to the second order differential equation

$$x'' - 2ax' + (a^2 + b^2)x = 0$$

together with $y = (ax - x')/b$. The general solution of (7) is easy to find:

$$x(t) = Ce^{at}\cos(bt - \gamma)\,x, \quad y(t) = Ce^{at}\sin(bt - \gamma) \tag{8}$$

where $C \geq 0$, $0 \leq \gamma < 2\pi$ are constants. This is a family of equiangular spirals depending on the two parameters $C$ and $\gamma$. When $a = 0$ we have concentric circles, which can be viewed as a particular case of (6) with $k = 0$ (Figure 4). When $b = 0$ we have rays going from the origin $O$ to infinity with slope $y/x = -\tan\gamma$ (Figure 5).

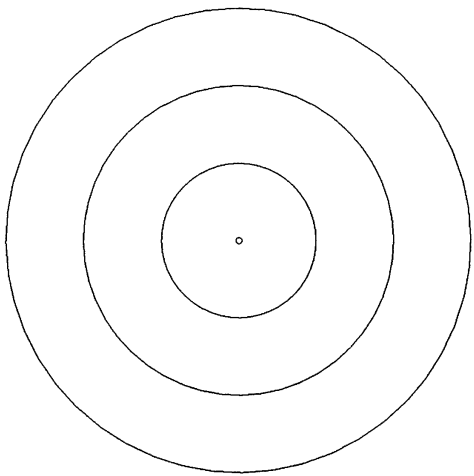
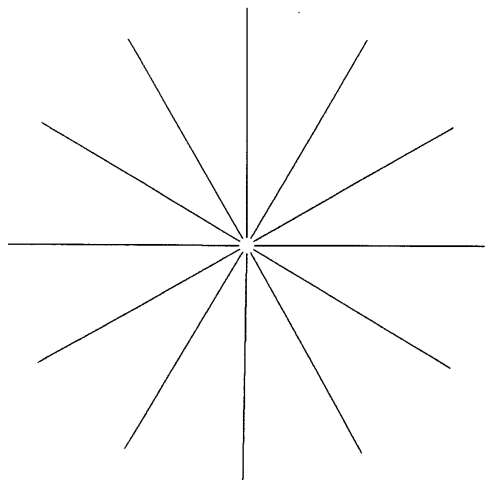

Figure 4                                                                          Figure 5

When $b \neq 0$, we can re-scale the parameter, $bt - \gamma \to t$ and setting $a/b = k$ we can write the solution in the form

$$x(t) = r_0 e^{kt} \cos t, \quad y(t) = r_0 e^{kt} \sin t \tag{9}$$

which is equivalent to equation (6) in polar coordinates, if we consider $t$ to be the polar angle. The proof is completed!

It is interesting to compute the circulation of the vector field (6) on some circle $G$ with center at the origin $O$. Keeping in mind that the function

$$\overline{V(z)} = u - iv = \frac{\overline{c}}{z} = \frac{a - ib}{z}$$

is holomorphic we write:

$$Cir(V:G) = \oint_G u\,dx + v\,dy = \text{Re} \oint_G (u - iv)\,d(x + iy) = \text{Re} \oint_G \frac{a - ib}{z}\,dz = 2\pi b$$

and the flux through that circle is

$$Flux(V:G) = \oint_G u\,dy - v\,dx = \text{Im} \oint_G (u - iv)\,d(x + iy) = \text{Im} \oint_G \frac{a - ib}{z}\,dz = 2\pi b$$

Therefore, $Cir(V:G)$ and $Flux(V:G)$ are independent of the radius of $G$ and also of each other. They are determined only by the complex constant $c = a + ib$ which can be arbitrary. When $a = 0$, the streamlines are concentric circles centered at $O$. We observe a vortex with no generation or disappearance of fluid. The case $a > 0$ corresponds to a source, and the case $a < 0$ to a sink. In our model, sources and sinks behave the same way, only the streamlines have different directions. It is clear from the equations (8) that when $a > 0$ the point $(x(t), y(t))$ moves away from the source at the origin $O$, and when $a < 0$ the point $(x(t), y(t))$ moves toward the sink at $O$. When $b = 0$ (no whirls!), the streamlines, as mentioned above, are radial rays; $b > 0$ corresponds to positive (counterclockwise) circulation (rotation about the vortex) and $b < 0$ indicates a clockwise circulation. The characteristic angle of the spirals is $\alpha = \text{arccot}(a/b)$, as represented in Figure 3.

## A simple tornado

Tornadoes and hurricanes are three dimensional formations. In order to construct a simple model of tornado, we extend the spiral motion from the horizontal $xy$-plane to three dimensions, along the vertical $z$-axis (*the letter $z$ plays now a different role - it will be used for the third Cartesian coordinate*). Our starting point is the system (7) which we replace by a $3 \times 3$ system of linear differential equations with real coefficients:

$$\begin{aligned} x' &= a_{11} x + a_{12} y + a_{13} z \\ y' &= a_{21} x + a_{22} y + a_{23} z \\ z' &= a_{31} x + a_{32} y + a_{33} z \end{aligned} \tag{10}$$

The characteristic equation of the real matrix $[a_{ij}]$ is a third order algebraic equation with real coefficients and it may have only real roots or two complex conjugate roots and one real. Since we want the $xy$-projection of the vector field $<x',y',z'>$ to satisfy two dimensional system (7), we exclude the case of all real roots, because this case brings to exponential growth or decay of the solutions in all dimensions. With two complex conjugate eigenvalues and one real eigenvalue, the system (10) splits into a direct product of one $2 \times 2$ real system of the form (7) for $x, y$ and one single equation for the third variable $z$, which can be put in the simple form .

$$z' = p\, z.$$

There are no other possibilities. The two complex eigenvalues bring equations (9) and the real eigenvalue, say, $p$, brings to the solution $z(t) = z_0 e^{pt}$ of the third equation. The general solution is

$$x(t) = r_0 e^{kt} \cos t$$
$$y(t) = r_0 e^{kt} \sin t$$
$$z(t) = z_0 e^{pt}$$

where $-\infty < t < \infty$, and $r_0$, $z_0$ are real constants. This is a family of spiral curves on the surface

$$x^2 + y^2 = M^2 z^{2k/p}$$

($M$ is a constant). When $k = p$ we have conchospirals, that is, spirals on the cone $x^2 + y^2 = M^2 z^2$ (see [2]). The tornado shape appears when we take $k < p < 0 < \min\{r_0, z_0\}$.

A streamline is shown in Figure 6 together with its projection on the $xy$-plane.

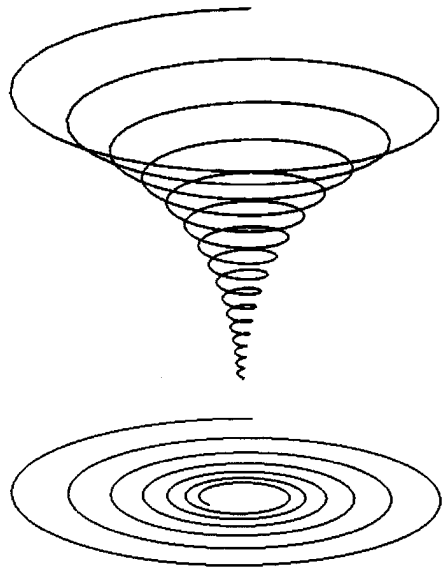

Figure 6

## Final notes

Bathtub vortices have been discussed, for instance, in [6]. Many books on complex variables include two-dimensional fluid flows. This is done very well in the classical book [10]. The topic is included sometimes in books on fluid mechanics [11, 12]. Vortices and sinks/sources, however, are often treated separately and the spiral shape does not appear ([10] and [11] are exceptions). A combination of one vortex and one sink, also called a spiral vortex, is mentioned as a model of tornado in [12, p. 231]: "A tornado may be approximated by a two-dimensional vortex and a sink, except in a region near the origin, and this is an example of good agreement between a real and ideal fluid".

In order to show that the vector field (5) has spiral streamlines we could integrate the function $\overline{V(z)} = \overline{c}/z$ and come to the complex potential $(a - ib)\log(z)$, which provides the combination of one vortex and a source/sink [3, 10]. Using the system (7) instead, we avoid the multivalued complex logarithm. Systems of differential equations have been used very efficiently to describe vector fields with one singularity ([1, 4]).

Real world tornados and hurricanes are very complex formations while our model provides only a rough approximation. More about tornados can be found in [9]. There is a resemblance to a spiral vortex in the form of some spiral galaxies like M51, M81, M100 and M101. University of Rochester researchers have fitted recently the spiral arms of M51 with pieces of equiangular spirals and discovered a good match in [7] (see Figure 7, courtesy of the authors).

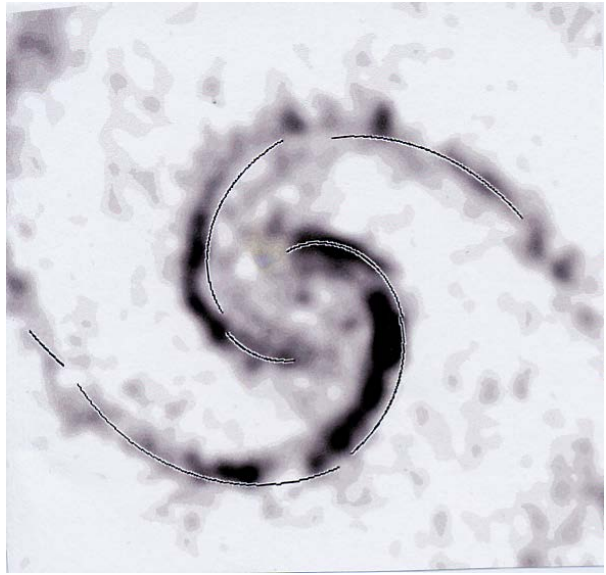

Figure 7. Galaxy M51

Some galaxies possibly behave like cosmic hurricanes in slow motion, where the density increases toward the center which plays the role of a sink/source. Interesting comments on this topic can be found on pp. 120-123 in [8].

*Acknowledgment*: Hurricane Bonnie and Galaxy M51 on Figures 1 and 2 are correspondingly from the NASA and the NASA/IPAC Extragalactic (NED) free databases.